\documentclass[10pt]{extarticle}

\usepackage[english]{babel}
\usepackage{graphicx}
\usepackage{framed}
\usepackage[normalem]{ulem}
\usepackage{indentfirst,ragged2e}
\usepackage{amsmath,amsthm,amssymb,amsfonts,wasysym,verbatim,bbm}
\usepackage{physics}
\usepackage[T1]{fontenc}
\usepackage{lmodern,mathrsfs}
\usepackage{pifont} 
\usepackage{mathdots} 
\usepackage{wrapfig}
\usepackage[inline,shortlabels]{enumitem}
\setlist{topsep=2pt,itemsep=2pt,parsep=0pt,partopsep=0pt}
\usepackage[dvipsnames]{xcolor}
\usepackage[utf8]{inputenc}
\usepackage[a4paper,top=0.5in,bottom=0.2in,left=0.5in,right=0.5in,footskip=0.3in,includefoot]{geometry} 
\usepackage{multicol}
\usepackage[most]{tcolorbox}
\usepackage{tikz,tikz-3dplot,tikz-cd,tkz-tab,tkz-euclide,tikzsymbols,pgf,pgfplots}
\pgfplotsset{compat=1.3}
\usepgfplotslibrary{fillbetween}
\usepackage[backend=bibtex,style=numeric]{biblatex}
\usepackage{hyperref}
\usepackage[nameinlink]{cleveref}

\definecolor{contcolor}{HTML}{514CE0}
\definecolor{convcolor}{HTML}{4CC6D4}

\newcommand{\ep}{\varepsilon}

\newcommand{\ceil}[1]{\ensuremath{\left\lceil#1\right\rceil}}

\newcommand{\C}{\mathbb{C}}

\newcommand{\N}{\mathbb{N}}
\newcommand{\Q}{\mathbb{Q}}
\newcommand{\R}{\mathbb{R}}

\newcommand{\Z}{\mathbb{Z}}
\newcommand{\As}{\mathcal{A}}

\newcommand{\Ms}{\mathcal{M}}

\newcommand{\limm}{\lim_{m\to\infty}}
\newcommand{\limn}{\lim_{n\to\infty}}

\newcommand{\summ}{\sum_{m=1}^\infty}

\newcommand{\exc}{\backslash}
\newcommand{\sub}{\subseteq}
\newcommand{\sups}{\supseteq}
\newcommand{\capp}{\bigcap}
\newcommand{\cupp}{\bigcup}

\newcommand{\cuppk}{\bigcup_{k=1}^\infty}

\newcommand{\cuppm}{\bigcup_{m=1}^\infty}

\newcommand{\measto}{\xrightarrow{\mu}}

\newtheoremstyle{mystyle}{}{}{}{}{\sffamily\bfseries}{.}{ }{}
\newtheoremstyle{cstyle}{}{}{}{}{\sffamily\bfseries}{.}{ }{\thmnote{#3}}

\theoremstyle{mystyle}{\newtheorem{definition}{Definition}
}
\theoremstyle{mystyle}{\newtheorem{proposition}[definition]{Proposition}}
\theoremstyle{mystyle}{\newtheorem{theorem}[definition]{Theorem}}
\theoremstyle{mystyle}{}
\theoremstyle{mystyle}{}
\theoremstyle{mystyle}{\newtheorem*{remark}{Remark}}
\theoremstyle{mystyle}{}
\theoremstyle{mystyle}{\newtheorem*{example}{Example}}
\theoremstyle{mystyle}{\newtheorem*{examples}{Examples}}
\theoremstyle{definition}{}
\theoremstyle{cstyle}{\newtheorem*{cthm}{}}
\newtheoremstyle{warn}{}{}{}{}{\normalfont}{}{ }{}
\theoremstyle{warn}

\newcommand{\warningsign}[1]{\tikz[scale=#1, every node/.style={scale=#1}]{\draw[-,line width={#1*0.8mm},red,fill=yellow,rounded corners={#1*2.5mm}] (0,0)--(1,{-sqrt(3)})--(-1,{-sqrt(3)})--cycle;
\node at (0,-1) {\fontsize{48}{60}\selectfont\bfseries!};}}

\newenvironment{talign*}{\csname align*\endcsname}{\endalign}

\usepackage[explicit]{titlesec}
\titleformat{\chapter}{\normalfont\sffamily\Huge\bfseries}{\thechapter}{20pt}{#1}
\titleformat{\section}{\normalfont\sffamily\Large\bfseries}{\thesection}{12pt}{#1}
\titleformat{\subsection}{\normalfont\sffamily\Large\bfseries}{\thesubsection}{12pt}{#1}
\titleformat{\subsubsection}{\normalfont\sffamily\large\bfseries}{\thesubsection}{8pt}{#1}

\titlespacing*{\chapter}{0pt}{-\topskip}{10pt}
\titlespacing*{\section}{0pt}{5pt}{5pt}
\titlespacing*{\subsection}{0pt}{5pt}{5pt}
\titlespacing*{\subsubsection}{0pt}{5pt}{5pt}

\DeclareMathAlphabet\mathbfcal{OMS}{cmsy}{b}{n}
\newcommand{\ind}{\hspace{0.2in}} 

\makeatletter
\g@addto@macro\normalsize{
\setlength\abovedisplayskip{3pt}
\setlength\belowdisplayskip{3pt}
\setlength\abovedisplayshortskip{0pt}
\setlength\belowdisplayshortskip{0pt}}
\makeatother

\makeatletter
\renewenvironment{abstract}{
\if@twocolumn
\section*{\abstractname}
\else
\begin{center}
{\sffamily\bfseries\abstractname\vspace{\z@}}
\end{center}
\quotation
\fi}
{\if@twocolumn\else\endquotation\fi}
\makeatother

\title{\Large\sffamily\bfseries A Result on Convergence of Double Sequences of Measurable Functions}
\author{\Large\sffamily Senan Sekhon}
\date{\large\sffamily April 20, 2021}

\begin{document}

\maketitle

\begin{abstract}
    We investigate a result on convergence of double sequences of numbers and how it extends to measurable functions.
\end{abstract}

\ind In this article, we will investigate the following statement:

\begin{cthm}[Statement 0]\label{statement}
Suppose $(a_{m,n})_{m,n=1}^\infty$ is a double sequence, $(b_m)_{m=1}^\infty$ is a sequence and $c$ is an element. Suppose $(a_{m,n})_{n=1}^\infty$ converges to $b_m$ as $n\to\infty$ for each $m\in\N$ and $(b_m)_{m=1}^\infty$ converges to $c$ as $m\to\infty$. Then there is a strictly increasing sequence $(n_m)_{m=1}^\infty$ in $\N$ such that $(a_{m,n_m})_{m=1}^\infty$ converges to $c$ as $m\to\infty$.
\end{cthm}

\begin{center}
\begin{tikzpicture}
\foreach \m in {1,2,3,4}
\foreach \n in {1,2,3,4,5,6}
\node (a\m\n) at (\n,-0.8*\m) {$a_{\m,\n}$};
\foreach \m in {1,2,3,4}
\node (ds\m) at (6.7,-0.8*\m) {$\cdots$};
\foreach \m in {1,2,3,4}
\node (b\m) at (8,-0.8*\m) {$b_\m$};
\node (dv) at (8,-3.7) {$\vdots$};
\node (c) at (8,-5) {$c$};
\foreach \m in {1,2,3,4}
\draw[->,red,shorten >=1mm] (ds\m)--(b\m);
\draw[->,blue,shorten <=1mm,shorten >=1mm] (dv)--(c);
\draw[->,thick] (0.8,-0.4)--(2,-0.4) node[midway,above]{$n$};
\draw[->,thick] (0.5,-0.8)--(0.5,-2) node[midway,left]{$m$};
\node (dg) at (5.5,-3.6) {$\ddots$};
\foreach \sx in {{(a11)--(a22)},{(a22)--(a34)},{(a34)--(a45)}}
\draw[->,Green] \sx;
\draw[->,Green,shorten >=2mm] (5.8,-4)--(c);
\node[text width=3.5cm] at (-2,-1.2) {If each \textit{row}\\\centering \textcolor{red}{$a_{m,1},a_{m,2},a_{m,3},...$}\\converges to some $b_m$,};
\node[text width=2.5cm] at (10,-1.2) {and the \textit{column}\\\centering \textcolor{blue}{$b_1,b_2,b_3,...$}\\converges to $c$,};
\node[text width=4.3cm] at (3.8,-4.8) {then there is some \textit{diagonal}\\\centering \textcolor{Green}{$a_{1,n_1},a_{2,n_2},a_{3,n_3},...$}\\that converges to $c$.};
\end{tikzpicture}
\end{center}

\ind Of course, the validity of this statement depends highly on the type of objects $a_{m,n}$, $b_m$ and $c$ are, as well as the notion of convergence in question. We first show that it holds in the case of complex numbers with the usual (Euclidean) notion of convergence:

\begin{proposition}\label{dsr}
Suppose $(a_{m,n})_{m,n=1}^\infty$ is a double sequence in $\C$, $(b_m)_{m=1}^\infty$ is a sequence in $\C$ and $c\in\C$. Suppose $a_{m,n}\to b_m$ as $n\to\infty$ for each $m\in\N$ and $b_m\to c$ as $m\to\infty$. Then there is a strictly increasing sequence $(n_m)_{m=1}^\infty$ in $\N$ such that $a_{m,n_m}\to c$ as $m\to\infty$.
\end{proposition}
\begin{proof}
Since $b_m\to c$ as $m\to\infty$, for any $\ep>0$, there exists $M_0\in\N$ such that $\abs{b_m-c}<\frac{\ep}{2}$ for all $m\ge M_0$. Now suppose $m\in\N$. Since $a_{m,n}\to b_m$ as $n\to\infty$, there exists $k_m\in\N$ such that $\abs{a_{m,n}-b_m}<\frac{1}{m}$ for all $n\ge k_m$. Define the sequence $(n_m)_{m=1}^\infty$ recursively by $n_1=k_1$ and $n_m=\max\{k_m,n_{m-1}+1\}$ for all $m\ge 2$. Then $(n_m)_{m=1}^\infty$ is a strictly increasing sequence in $\N$, and $\abs{a_{m,n_m}-b_m}<\frac{1}{m}$ for all $m\in\N$. Now define $M=\max\left\{M_0,\ceil{\frac{2}{\ep}}+1\right\}$. Then for all $m\ge M$, we have $\abs{a_{m,n_m}-c}\le\abs{a_{m,n_m}-b_m}+\abs{b_m-c}<\frac{\ep}{2}+\frac{\ep}{2}=\ep$. Thus $a_{m,n_m}\to c$ as $m\to\infty$.
\end{proof}

The following examples show that the ``straight'' diagonal subsequence $(a_{m,m})_{m=1}^\infty$ may not converge to $c$:
\begin{examples}\leavevmode
\begin{enumerate}
    \item Suppose $a_{m,n}=\frac{m}{n}$. Then for each $n\in\N$, we have $\limn a_{m,n}=\limn\frac{m}{n}=0$, so $b_m=0$ for all $m\in\N$. Thus $\limm b_m=\limm 0=0$, so $c=0$. However, since $a_{m,m}=\frac{m}{m}=1$ for all $m\in\N$, we have $\limm a_{m,m}=1\ne0$. Instead, we can set $n_m=m^2$ to get $a_{m,n_m}=\frac{m}{m^2}=\frac{1}{m}\to0$ as $m\to\infty$.
    \item Suppose $a_{m,n}=\frac{m^2}{n}$. Again, we have $b_m=0$ for all $m\in\N$, and $c=0$. However, since $a_{m,m}=\frac{m^2}{m}=m$ for all $m\in\N$, we have $\limm a_{m,m}=\infty$. Instead, we can set $n_m=m^3$ to get $a_{m,n_m}=\frac{m^2}{m^3}=\frac{1}{m}\to0$ as $m\to\infty$.
\end{enumerate}
\end{examples}
Additionally, the last example above shows that despite the assumptions of \Cref{dsr}, the double sequence $(a_{m,n})_{m,n=1}^\infty$ may be unbounded.

\begin{remark}
If $a_{m,n}\to b_m$ \textit{uniformly} in $m$, i.e. for any $\ep>0$, there exists $N\in\N$ (independent of $m$) such that $\abs{a_{m,n}-b_m}<\ep$ for all $n\ge N$ and all $m\in\N$, then there is a straight diagonal subsequence $(a_{m,m+p})_{m=1}^\infty$ that converges to $c$. This can be seen from the proof: If $a_{m,n}\to b_m$ uniformly in $m$, then $k_m$ will not depend on $n$, and the sequence $(n_m)_{m=1}^\infty$ will simply be $(k_1,k_1+1,k_1+2,k_1+3,...)$, so we can set $p=k_1$. However, we may still need to `translate' the diagonal, i.e. use $(a_{m,m+p})_{m=1}^\infty$ instead of $(a_{m,m})_{m=1}^\infty$.
\end{remark}

\subsection*{Metric and Topological Convergence}
We now show that \hyperref[statement]{Statement 0} holds for convergence in any metric space:

\begin{theorem}\label{dsmet}
Suppose $X$ is a metric space, $(x_{m,n})_{m,n=1}^\infty$ is a double sequence in $X$, $(y_m)_{m=1}^\infty$ is a sequence in $X$ and $z\in X$. Suppose $x_{m,n}\to y_m$ as $n\to\infty$ for each $m\in\N$ and $y_m\to z$ as $m\to\infty$. Then there is a strictly increasing sequence $(n_m)_{m=1}^\infty$ in $\N$ such that $x_{m,n_m}\to z$ as $m\to\infty$.
\end{theorem}
\begin{proof}
Since $y_m\to z$ as $m\to\infty$, for any $\ep>0$, there exists $M_0\in\N$ such that $d(y_m,z)<\frac{\ep}{2}$ for all $m\ge M_0$. Now suppose $m\in\N$. Since $x_{m,n}\to y_m$ as $n\to\infty$, there exists $k_m\in\N$ such that $d(x_{m,n},y_m)<\frac{1}{m}$ for all $n\ge k_m$. Define the sequence $(n_m)_{m=1}^\infty$ recursively by $n_1=k_1$ and $n_m=\max\{k_m,n_{m-1}+1\}$ for all $m\ge 2$. Then $(n_m)_{m=1}^\infty$ is a strictly increasing sequence in $\N$, and $d(x_{m,n_m},y_m)<\frac{1}{m}$ for all $m\in\N$. Now define $M=\max\left\{M_0,\ceil{\frac{2}{\ep}}+1\right\}$. Then for all $m\ge M$, we have $d(x_{m,n_m},z)\le d(x_{m,n_m},y_m)+d(y_m,z)<\frac{\ep}{2}+\frac{\ep}{2}=\ep$. Thus $x_{m,n_m}\to z$ as $m\to\infty$.
\end{proof}

\ind In particular, \hyperref[statement]{Statement 0} holds for \emph{uniform} convergence of sequences of complex-valued functions $f:S\to\C$, since this is equivalent to convergence in the supremum metric: $d_\infty(f,g)=\sup_{x\in S} \abs{f(x)-g(x)}$.\\

\ind The natural question now is whether \hyperref[statement]{Statement 0} holds for \emph{pointwise} convergence. We now present some examples to show that it does not.\\

\ind The first example is from \cite[Exercise 1.5.7, Page 125]{tao}.
\begin{example}
Suppose $X=\N^\N$ (the set of all sequences of positive integers). Define $h=0$ and $g_m=0$ for all $m\in\N$. Then $g_m\to h$ pointwise as $m\to\infty$. Now for each $m,n\in\N$, define $f_{m,n}=1_{A_{m,n}}$, where $A_{m,n}$ is the set of all sequences $(a_s)_{s=1}^\infty$ in $\N^\N$ such that $a_m\ge n$. Then $f_{m,n}\to0$ pointwise as $n\to\infty$ for each $m\in\N$. Now suppose $(n_m)_{m=1}^\infty$ is a strictly increasing sequence in $\N$. Then $f_{m,n_m}((n_s)_{s=1}^\infty)=1$ for all $m\in\N$ (since $n_m\ge n_m$). Thus $f_{m,n_m}((n_s)_{s=1}^\infty)\nrightarrow0$ pointwise as $m\to\infty$.
\end{example}

\ind The next example is slightly more involved.
\begin{example}
Suppose $(q_m)_{m=1}^\infty$ is an enumeration of $\Q$. Define $h=0$, and for each $m\in\N$, define $g_m=1_{\{q_m\}}$. Then $g_m\to h$ pointwise as $m\to\infty$. Now for each $m,n\in\N$, define $f_{m,n}=1_{U_{m,n}}$, where $U_{m,n}=\left(q_m-\frac{1}{n},q_m+\frac{1}{n}\right)$. Then $f_{m,n}\to g_m$ pointwise as $n\to\infty$ for each $m\in\N$. Now suppose $(n_m)_{m=1}^\infty$ is a strictly increasing sequence in $\N$. We want to show that $f_{m,n_m}\nrightarrow0$ pointwise as $m\to\infty$. We will do this by showing that there exists $r\in\R$ such that $f_{m,n_m}(r)=1$ for infinitely many $m\in\N$.\\

Since $U_{1,n_1}$ is an open interval, it contains some closed interval $[a_1,b_1]$. This interval contains a rational number $q_{k_2}\ne q_{k_1}$, and the corresponding open interval $U_{k_2,n_{k_2}}$ contains a closed interval $[a_2,b_2]$. Assume without loss of generality that $[a_2,b_2]\sub[a_1,b_1]$ (otherwise, replace $[a_2,b_2]$ with $[a_1,b_1]\cap[a_2,b_2]$). This interval contains a rational number $q_{k_3}\ne q_{k_1},q_{k_2}$, and the corresponding open interval $U_{k_3,n_{k_3}}$ contains a closed interval $[a_3,b_3]$ (which we can assume is contained in $[a_2,b_2]$). Repeating this process yields a decreasing sequence $([a_l,b_l])_{l=1}^\infty$ of closed intervals such that $[a_l,b_l]\sub U_{k_l,n_{k_l}}$ for each $l\in\N$. Note that all $k_l$ are distinct by construction, since $q_{k_l}\ne q_{k_1},q_{k_2},...,q_{k_{l-1}}$. Since $([a_l,b_l])_{l=1}^\infty$ is a decreasing sequence of closed bounded intervals, their intersection is non-empty, i.e. there exists $r\in\capp_{l=1}^\infty [a_l,b_l]$. Thus $r\in[a_l,b_l]\sub U_{k_l,n_{k_l}}$ for all $l\in\N$. In particular, $r\in U_{m,n_m}$ for infinitely many $m\in\N$, and so $f_{m,n_m}(r)=1$ for infinitely many $m\in\N$.\\

Thus $f_{m,n_m}(r)\nrightarrow0$ as $m\to\infty$, and so $f_{m,n_m}\nrightarrow0$ pointwise as $m\to\infty$.
\end{example}

\ind The next example uses the fact that the set of discontinuities of any Baire class 1 function is meager in $\R$ (and thus not equal to $\R$, by the Baire category theorem). See \cite{fung} for details.
\begin{example}
Define $h=1_\Q$, and for each $m\in\N$, define $g_m=1_{A_m}$, where $A_m=\frac{1}{m!}\Z=\left\{\frac{k}{m!}\mid k\in\Z\right\}$. Then $g_m\to h$ pointwise as $m\to\infty$. Now for each $m,n\in\N$, define $f_{m,n}(x)=\cos(m!\pi x)^{2n}$. Then $f_{m,n}\to g_m$ pointwise as $n\to\infty$ for each $m\in\N$. Note that all $f_{m,n}$ are continuous everywhere on $\R$, but $h$ is discontinuous everywhere on $\R$, so it is of Baire class 2. Thus there cannot be a diagonal sequence $f_{m,n_m}$ that converges pointwise to $h$.
\end{example}

\ind We now extend \Cref{dsmet} to first countable topological spaces. Note that while limits of sequences in such spaces may not be unique, \hyperref[statement]{Statement 0} still holds (and it holds for all of the limits in question).

\begin{theorem}\label{dsfc}
Suppose $X$ is a first countable topological space, $(x_{m,n})_{m,n=1}^\infty$ is a double sequence in $X$, $(y_m)_{m=1}^\infty$ is a sequence in $X$ and $z\in X$. Suppose $x_{m,n}\to y_m$ as $n\to\infty$ for each $m\in\N$ and $y_m\to z$ as $m\to\infty$. Then there is a strictly increasing sequence $(n_m)_{m=1}^\infty$ in $\N$ such that $x_{m,n_m}\to z$ as $m\to\infty$.
\end{theorem}
\begin{proof}
Suppose $(U_k)_{k=1}^\infty$ is a local base at $z$. Assume without loss of generality that $U_k\sups U_{k+1}$ for all $k\in\N$ (otherwise, replace each $U_k$ with $\capp_{p=1}^k U_p$). Since $y_m\to z$ as $m\to\infty$, for any $k\in\N$, there exists $M_k\in\N$ such that $y_m\in U_k$ for all $m\ge M_k$. Now suppose $m\in\N$ and $(V_{m,r})_{r=1}^\infty$ is a local base at $y_m$. Again, assume without loss of generality that $V_{m,r}\sups V_{m,r+1}$ for all $m,r\in\N$. Since $x_{m,n}\to y_m$ as $n\to\infty$, for any $r\in\N$, there exists $N_{m,r}\in\N$ such that $x_{m,n}\in V_{m,r}$ for all $n\ge N_{m,r}$. Define the sequence $(n_s)_{s=1}^\infty$ recursively by $n_1=N_{1,1}$ and $n_s=\max\left\{N_{M_s,s},n_{s-1}\right\}$ for all $s\ge2$.\\

Now suppose $W$ is an open neighborhood of $z$. Since $(U_k)_{k=1}^\infty$ is a local base at $z$, there exists $K_0\in\N$ such that $U_{K_0}\sub W$. This implies that $y_m\in W$ for all $m\ge M_{K_0}$. Now suppose $m\in\N$, $m\ge M_{K_0}$. Since $(V_{m,r})_{r=1}^\infty$ is a local base at $y_m$, there exists $L_{0,m}\in\N$ such that $V_{m,L_{0,m}}\sub W$. This implies that $x_{m,n}\in W$ for all $m\ge M_{K_0}$ and all $n\ge N_{M_{K_0},L_{0,M_{K_0}}}$. Now define $S_0$ as the least $s\in\N$ such that $s\ge K_0$ and $N_{M_s,s}\ge N_{M_{K_0},L_{0,M_{K_0}}}$. Then for all $s\ge S_0$, we have $x_{s,n_s}\in W$. Thus $x_{s,n_s}\to z$ as $s\to\infty$.
\end{proof}

\ind While \hyperref[statement]{Statement 0} does not hold for pointwise convergence in general, it does if the domain is a countable set. This is because if $S$ is a countable set and $X$ is a first countable topological space, then the set of functions $f:S\to X$ with pointwise convergence (or equivalently, the set $X^S$ with the product topology) is also first countable. See \cite[Theorem 7.1.7, Page 172]{singh} for a proof.

\subsection*{Measure-Theoretic Convergence}
We now turn to measure-theoretic notions of convergence. From now on, $\Ms(X)$ denotes the set of all complex-valued measurable functions $f:X\to\C$, where $(X,\As,\mu)$ is an arbitrary measure space (not assumed to be $\sigma$-finite).

\begin{definition}
Suppose $(X,\As,\mu)$ is a measure space, $(f_n)_{n=1}^\infty$ is a sequence of functions in $\Ms(X)$ and $f\in\Ms(X)$.
\begin{itemize}
    \item $(f_n)$ \textbf{converges almost uniformly} to $f$ (denoted by $f_n\to f$ a.u.) if for any $\ep>0$, there is a set $A\in\As$ such that $\mu(A)<\ep$ and $f_n\to f$ uniformly on $X\exc A$.
    \item $(f_n)$ \textbf{converges almost everywhere} to $f$ (denoted by $f_n\to f$ a.e.) if there is a set $N\in\As$ such that $\mu(N)=0$ and $f_n\to f$ pointwise on $X\exc N$.
    \item $(f_n)$ \textbf{converges in measure} to $f$ (denoted by $f_n\measto f$) if for any $\ep>0$, we have $\limn\mu(\{x\in X\mid\abs{f_n(x)-f(x)}\ge\ep\})=0$.
\end{itemize}
\end{definition}

\ind \hyperref[statement]{Statement 0} automatically holds for convergence in measure, as it is a metric convergence (after identifying functions that are equal a.e.), given by the \emph{Fréchet metric}: $d(f,g)=\inf_{r>0} (\mu(\{\abs{f-g}\ge r\})+r)$. Similarly, it holds for (strong) convergence in $L^p$ for any $1\le p\le\infty$, as this is a norm (and thus metric) convergence.\\

\ind We now show that \hyperref[statement]{Statement 0} holds for a.u. convergence:

\begin{theorem}\label{dsau}
Suppose $(X,\As,\mu)$ is a measure space, $(f_{m,n})_{m,n=1}^\infty$ is a double sequence in $\Ms(X)$, $(g_m)_{m=1}^\infty$ is a sequence in $\Ms(X)$ and $h\in\Ms(X)$. Suppose $f_{m,n}\to g_m$ a.u. as $n\to\infty$ for each $m\in\N$ and $g_m\to h$ a.u. as $m\to\infty$. Then there is a strictly increasing sequence $(n_m)_{m=1}^\infty$ in $\N$ such that $f_{m,n_m}\to h$ a.u. as $m\to\infty$.
\end{theorem}
\begin{proof}
Since $g_m\to h$ a.u. as $m\to\infty$, for any $\delta>0$, there is a set $A\in\As$ such that $\mu(A)<\frac{\delta}{2}$ and $g_m\to h$ uniformly on $X\exc A$ as $m\to\infty$. Now suppose $m\in\N$. Since $f_{m,n}\to g_m$ a.u., there is a set $B_m\in\As$ such that $\mu(B_m)<\frac{\delta}{2^{m+1}}$ and $f_{m,n}\to g_m$ uniformly on $X\exc B_m$ as $n\to\infty$. Define $C=A\cup(\cuppm B_m)$. Then $\mu(C)\le\mu(A)+\summ\mu(B_m)<\frac{\delta}{2}+\frac{\delta}{2^{m+1}}=\frac{\delta}{2}+\frac{\delta}{2}=\delta$ and on $X\exc C$, we have $f_{m,n}\to g_m$ uniformly for all $m\in\N$ and $g_m\to h$ uniformly. By \hyperref[statement]{Statement 0} (for uniform convergence), there is a strictly increasing sequence $(n_m)_{m=1}^\infty$ in $\N$ such that $f_{m,n_m}\to h$ uniformly on $X\exc C$ as $m\to\infty$. At this point, we cannot conclude that $f_{m,n_m}\to h$ a.u., since the sequence $(n_m)_{m=1}^\infty$ may depend on $\delta$.\\

By what we have just shown, for each $k\in\N$, there is a set $C_k\in\As$ and a strictly increasing sequence $(n_{k,m})_{m=1}^\infty$ in $\N$ such that $\mu(C_k)<2^{-k}$ and $f_{m,n_{k,m}}\to h$ uniformly on $X\exc C_k$ as $m\to\infty$. For each $k\in\N$, define $D_k=\cupp_{l=k}^\infty C_l$. Then $(D_k)_{k=1}^\infty$ is a decreasing sequence of sets in $\As$ and $\mu(D_k)\le\sum_{l=k}^\infty \mu(C_l)<\sum_{l=k}^\infty 2^{-l}=2^{1-k}$. As we have already shown, for each $k\in\N$, there is a strictly increasing sequence $(n_{k,m})_{m=1}^\infty$ in $\N$ such that $f_{m,n_{k,m}}\to h$ uniformly on $X\exc C_k$ (and thus on $X\exc D_k$, since $D_k\sups C_k$).\\

As we have just shown, there is a set $D_1\in\As$ and a strictly increasing sequence $(n_{1,m})_{m=1}^\infty$ in $\N$ such that $\mu(D_1)<1$ and $f_{m,n_{1,m}}\to h$ uniformly on $X\exc D_1$ as $m\to\infty$. For each $m\in\N$, since $(f_{m,n_{1,p}})_{p=1}^\infty$ is a subsequence of $(f_{m,n})_{n=1}^\infty$, we have $f_{m,n_{1,p}}\to g_m$ a.u. as $p\to\infty$. In other words, the double sequence $(f_{m,n_{1,p}})_{m,p=1}^\infty$ satisfies the same assumptions as $(f_{m,n})_{m,n=1}^\infty$. Thus there is a set $D_2\in\As$ and a subsequence $(n_{2,m})_{m=1}^\infty$ of $(n_{1,m})_{m=1}^\infty$ such that $\mu(D_2)<\frac{1}{2}$ and $f_{m,n_{2,m}}\to h$ uniformly on $X\exc D_2$ as $m\to\infty$. Repeating this process yields a decreasing sequence $(D_k)_{k=1}^\infty$ of sets in $\As$ and a nested sequence $((n_{k,m})_{m=1}^\infty)_{k=1}^\infty$ of strictly increasing sequences in $\N$ (i.e. for each $k\in\N$, $(n_{k+1,m})_{m=1}^\infty$ is a subsequence of $(n_{k,m})_{m=1}^\infty$) such that $\mu(D_k)<2^{1-k}$ and $f_{m,n_{k,m}}\to h$ uniformly on $X\exc D_k$ as $m\to\infty$. Now consider the diagonal sequence $(n_{m,m})_{m=1}^\infty$. For each $k\in\N$, this is a subsequence of $(n_{k,m})_{m=1}^\infty$ (except possibly for the first $k-1$ terms, which does not affect convergence), and as we have shown, $f_{m,n_{k,m}}\to h$ uniformly on $X\exc D_k$ as $m\to\infty$. Thus for all $k\in\N$, we have $f_{m,n_{m,m}}\to h$ uniformly on $X\exc D_k$ as $m\to\infty$.\\

Finally, suppose $\ep>0$ and set $k=\ceil{\log_2\left(\frac{1}{\ep}\right)}+2$. By what we have just shown, we have $\mu(D_k)<2^{1-k}<\ep$ and $f_{m,n_{m,m}}\to h$ uniformly on $X\exc D_k$. Thus $f_{m,n_{m,m}}\to f$ a.u.
\end{proof}

\ind What about a.e. convergence? Clearly \hyperref[statement]{Statement 0} cannot hold for a.e. convergence in general, since we can use any of the counterexamples for pointwise convergence (with the counting measure, a.e. convergence is equivalent to pointwise convergence). However, as we pointed out earlier, \hyperref[statement]{Statement 0} can be recovered for pointwise convergence if we assume the domain is a countable set. This suggests that we might be also be able to recover it for a.e. convergence if we assume the domain is $\sigma$-finite. We now show that this is indeed true:

\begin{theorem}
Suppose $(X,\As,\mu)$ is a $\sigma$-finite measure space, $(f_{m,n})_{m,n=1}^\infty$ is a double sequence in $\Ms(X)$, $(g_m)_{m=1}^\infty$ is a sequence in $\Ms(X)$ and $h\in\Ms(X)$. Suppose $f_{m,n}\to g_m$ a.e. as $n\to\infty$ for each $m\in\N$ and $g_m\to h$ a.e. as $m\to\infty$. Then there is a strictly increasing sequence $(n_m)_{m=1}^\infty$ in $\N$ such that $f_{m,n_m}\to h$ a.e. as $m\to\infty$.
\end{theorem}
\begin{proof}
Since $\mu$ is $\sigma$-finite, there is an increasing sequence $(X_k)_{k=1}^\infty$ of sets in $\As$ such that $\mu(X_k)<\infty$ for all $k\in\N$ and $\cuppk X_k=X$. For each $k\in\N$, by Egorov's theorem, we have $f_{m,n}\to g_m$ a.u. on $X_k$ as $n\to\infty$ for each $m\in\N$ and $g_m\to h$ a.u. on $X_k$ as $m\to\infty$. By \Cref{dsau}, there is a strictly increasing sequence $(n_{1,m})_{m=1}^\infty$ in $\N$ such that $f_{m,n_{1,m}}\to h$ a.u. on $X_1$ as $m\to\infty$. For each $m\in\N$, since $(f_{m,n_{1,p}})_{p=1}^\infty$ is a subsequence of $(f_{m,n})_{n=1}^\infty$, we have $f_{m,n_{1,p}}\to g_m$ a.u. as $p\to\infty$. In other words, the double sequence $(f_{m,n_{1,p}})_{m,p=1}^\infty$ satisfies the same assumptions as $(f_{m,n})_{m,n=1}^\infty$. Thus there is a subsequence $(n_{2,m})_{m=1}^\infty$ of $(n_{1,m})_{m=1}^\infty$ such that $f_{m,n_{2,m}}\to h$ a.u. as $m\to\infty$. Repeating this process yields a nested sequence $((n_{k,m})_{m=1}^\infty)_{k=1}^\infty$ of strictly increasing sequences in $\N$ (i.e. for each $k\in\N$, $(n_{k+1,m})_{m=1}^\infty$ is a subsequence of $(n_{k,m})_{m=1}^\infty$) such that $f_{m,n_{k,m}}\to h$ a.u. on $X_k$. Now consider the diagonal sequence $(n_{m,m})_{m=1}^\infty$. For each $k\in\N$, this is a subsequence of $(n_{k,m})_{m=1}^\infty$ (except possibly for the first $k-1$ terms, which does not affect convergence), and as we have shown, $f_{m,n_{k,m}}\to h$ a.u. on $X_k$ as $m\to\infty$. Thus for all $k\in\N$, we have $f_{m,n_{m,m}}\to h$ a.u. on $X_k$ as $m\to\infty$. Finally, since $f_{m,n_{m,m}}\to h$ a.u. (and thus a.e.) on $X_k$ as $m\to\infty$, we have $f_{m,n_{m,m}}\to h$ a.e. on $X$ as $m\to\infty$.
\end{proof}

\subsection*{Other Types of Convergence}
Finally, we turn to convergence in normed vector spaces. \hyperref[statement]{Statement 0} automatically holds for (strong) convergence in any normed vector space, since this is a metric convergence.\\

\ind We now give an example to show that \hyperref[statement]{Statement 0} does not hold for weak convergence in a normed vector space (not even in a Hilbert space):

\begin{example}
Suppose $X=\ell^2$ (the set of all square-summable sequences of real numbers) and $(e_k)_{k=1}^\infty$ is the standard basis of $\ell^2$, i.e. $e_k$ is the sequence with $1$ in the $k$\textsuperscript{th} place and $0$ elsewhere. Define $z=0$ and $y_m=0$ for all $m\in\N$. Then $y_m\to0$ strongly (and thus weakly) as $m\to\infty$. Now for each $m,n\in\N$, define $x_{m,n}=me_n$. Then for all $m\in\N$, we have $x_{m,n}\to0$ weakly as $n\to\infty$ (since $(me_n)_{n=1}^\infty$ is an orthonormal sequence scaled by $m$). However, for any sequence $(n_m)_{m=1}^\infty$ in $\N$, the corrsponding sequence $(x_{m,n_m})_{m=1}^\infty$ is unbounded, and so it cannot converge weakly.
\end{example}
\begin{remark}
This counterexample also works for weak convergence in $\ell^p$ for any $1\le p<\infty$, thereby showing that \hyperref[statement]{Statement 0} does not hold for weak convergence in $L^p$.
\end{remark}

\ind Note that the previous example was possible as the double sequence $(x_{m,n})_{m,n=1}^\infty$ was unbounded.\\

\ind If $X$ is a separable Banach space, then the weak* topology on $X^*$ is metrizable on bounded sets. If $X^*$ is separable, then the weak topology on $X$ is also metrizable on bounded sets (see \cite[Chapter V, Theorem 5.1, Pages 134-135]{conway} for details). Thus \hyperref[statement]{Statement 0} holds for weak convergence (if $X^*$ is separable and $(x_{m,n})_{m,n=1}^\infty$ is bounded) and for weak* convergence (if $X$ is separable and $(x_{m,n})_{m,n=1}^\infty$ is bounded).\\

\ind \hyperref[statement]{Statement 0} also holds for weak convergence in $\ell^1$, even if $(x_{m,n})_{m,n=1}^\infty$ is unbounded. This is because for sequences (but not nets) in $\ell^1$, weak convergence is equivalent to strong convergence. See \cite[Chapter V, Proposition 5.2, Pages 135-136]{conway} for a proof.

\vspace{5mm}

\subsubsection*{Acknowledgments}
I would like to thank Alex Derkach and James Baxter for their contributions.

\vfill

\printbibliography

\end{document}